\documentclass [12pt] {article}
\usepackage {amsmath}
\usepackage {amsfonts, amssymb,amsthm}
\usepackage {latexsym}
\usepackage[cp866]{inputenc}
\usepackage[russian]{babel}
\usepackage {amssymb}
\usepackage {amsmath}
\usepackage {graphicx}
\usepackage {longtable}

\begin{document}
\vskip1mm

\bigskip

\begin{center}
{\bf  ¡ ®¤­®¬ ªà¨â¥à¨¨ á®¡áâ¢¥­­®áâ¨ ¬­®£®§­ ç­ëå ®â®¡à ¦¥­¨©}
\end{center}

\begin{center}
{\bf ®à®â­¨ª®¢ ..\footnote{ ¡®â  ¯®¤¤¥à¦ ­  £à ­â®¬ 
ü01-01-00425}\footnote{Current address: CMUC, Apartado 3008, 3001 - 454 Coimbra, Portugal, mitvorot@mat.uc.pt}}
\end{center}

§¢¥áâ­  â¥á­ ï á¢ï§ì ¬¥¦¤ã á®¡áâ¢¥­­®áâìî ¨ § ¬ª­ãâ®áâìî
®¤­®§­ ç­ëå ®â®¡à ¦¥­¨©.  ¯à¨¬¥à [1], ¤«ï ®â®¡à ¦¥­¨© ¡ ­ å®¢ëå
¡¥áª®­¥ç­®¬¥à­ëå ¬­®£®®¡à §¨©, ®¡« ¤ îé¨å á¢®©áâ\-¢®¬
äà¥¤£®«ì¬®¢®áâ¨, íâ¨ ¯®­ïâ¨ï íª¢¨¢ «¥­â­ë.  â® ¦¥ ¢à¥¬ï á¢®©áâ¢®
á®¡áâ¢¥­­®áâ¨ ®â®¡à \-¦¥­¨ï ï¢«ï¥âáï ®¤­¨¬ ¨§ æ¥­âà «ì­ëå ¯à¨
¯®áâà®¥­¨¨ à §«¨ç­ëå â¥®à¨© â®¯®«®£¨ç¥áª®© áâ¥¯¥­¨.

 íâ®© à ¡®â¥ ¯à¨¢®¤¨âáï ®¡é¨© ªà¨â¥à¨© á®¡áâ¢¥­­®áâ¨ ¤«ï ¬­®£®§­ ç­ëå
®â®¡à ¦¥­¨©, ¡«¨§ª¨© ¯® ¤ãåã ª ¢ëè¥ã¯®¬ï­ãâ®¬ã á¢®©áâ¢ã äà¥¤£®«ì¬®¢ëå
®â®¡à ¦¥­¨©. à¨â¥à¨© ¯à¨¢¥¤¥­ ¢ ä®à¬¥ ¤¢ãå â¥®à¥¬: ¯àï¬®© ¨ ®¡à â­®©.

\section{ ¬ã«ìâ¨®â®¡à ¦¥­¨ïå}

 á¢®©áâ¢ å ¬ã«ìâ¨®â®¡à ¦¥­¨© ¬®¦­® ¯®¤à®¡­® ã§­ âì, ­ ¯à¨¬¥à, ¢ [2]. ë ¦¥
§¤¥áì ®â¬¥â¨¬ «¨èì í«¥¬¥­â à­ë¥ ¢®¯à®áë.

ãáâì $X,Y$ - ¬¥âà¨ç¥áª¨¥ ¯à®áâà ­áâ¢ . ­®£®§­ ç­ë¬ ®â®¡à ¦¥­¨¥¬
$F$ ¨§ $X$ ¢ $Y$ (®¡®§­ ç \-¥âáï $F:X \multimap Y$) ­ §ë¢ ¥âáï
á®®â¢¥âáâ¢¨¥, á®¯®áâ ¢«ïîé¥¥ ª ¦¤®© â®çª¥ $x \in X$ ­¥ª®¥ ­¥¯ãáâ®¥
¯®¤¬­®¦¥áâ¢® ¯à®áâà ­áâ¢  $Y$, ®¡®§­ ç ¥¬®¥ $F(x)$. ­®£®§­ ç­ë¥
®â®¡à ¦¥­¨ï ­ §ë¢ îâ â ª¦¥ ¬ã«ìâ¨®â®¡à ¦¥­¨ï¬¨.

«ï $A \subset X${\it ®¡à §®¬} $F(A)$ ­ §ë¢ ¥âáï ¬­®¦¥áâ¢® $\{ y
\in Y|\exists x \in A,$ $F(x) \ni y\} $.

«ï $B \subset Y${\it ¯®«­ë¬ ¯à®®¡à §®¬} $F_{ - }^{ - 1} (B)$
­ §ë¢ ¥âáï ¬­®¦¥áâ¢® $\{ x \in X|F(x) \cap B \ne \emptyset \} $.

 ç áâ­®áâ¨, ¤«ï $y\in Y:F_{-}^{-1}(y)=\{x\in X|y\in F(x)\}$.

á«¨ ¤«ï «î¡®£® § ¬ª­ãâ®£® $B \subset Y$ ¬­®¦¥áâ¢® $F_{ - }^{ - 1}
(B)$ § ¬ª­ãâ®, â® $F$ ­ §ë¢ ¥âáï {\it ¯®«ã­¥¯à¥\-àë¢­ë¬ á¢¥àåã}.

á«¨ $F$ ¯®«ã­¥¯à¥àë¢­® á¢¥àåã ¨ ®¡à § $F(x)$ ¢áïª®© â®çª¨ $x \in
X$ ª®¬¯ ªâ¥­, â® ®¡à § $F(A)$ ¢áïª®£® ª®¬¯ ªâ­®£® $A \subset X$
ª®¬¯ ªâ¥­ [2].

á«¨ ¤«ï «î¡®£® ®âªàëâ®£® $B \subset Y$ ¬­®¦¥áâ¢® $F_{ - }^{ - 1}
(B)$ ®âªàëâ®, â® $F$ ­ §ë¢ ¥âáï {\it ¯®«ã­¥¯à¥\-àë¢­ë¬ á­¨§ã}.

â®¡à ¦¥­¨¥ $F$ ­ §ë¢ ¥âáï {\it § ¬ª­ãâë¬}, ¥á«¨ ¤«ï «î¡ëå
¯®á«¥¤®¢ â¥«ì­®áâ¥© $x_{n} \mathop { \to }\limits_{n \to \infty }
x_{0} $ ¨§ $X$ ¨ $y_n \mathop { \to }\limits_{n \to \infty } y_{0}
$ ¨§ $Y$, $y_{n} \in F(x_{n} )$ ¢ë¯®«­¥­® $y_{0}^{} \in F(x_{0}^{}
)$.

áïª®¥ ¯®«ã­¥¯à¥àë¢­®¥ á¢¥àåã ¬ã«ìâ¨®â®¡à ¦¥­¨¥, ã ª®â®à®£® ®¡à § ª ¦¤®©
â®çª¨ § ¬ª­ãâ, ï¢«ï¥âáï § ¬ª­ãâë¬ [2].

ã¤¥¬ ­ §ë¢ âì ¬­®£®§­ ç­®¥ ®â®¡à ¦¥­¨¥ $F:X \multimap Y${\it
á®¡áâ¢¥­\-­ë¬}, ¥á«¨ ¤«ï «î¡®£® ª®¬¯ ªâ­®£® $B \subset Y$
¬­®¦¥áâ¢® $F_{-}^{-1}(B)$ ª®¬¯ ªâ­®.

ã¤¥¬ ­ §ë¢ âì ®â®¡à ¦¥­¨¥ $F$ {\it â®¯®«®£¨ç¥áª¨ § ¬ª­ãâë¬}, ¥á«¨
¤«ï «î¡®£® § ¬ª­ãâ®£® $A \subset X$ § ¬ª­ãâë¬ ¡ã¤¥â ¨ $F(A)$.

{\bf  ¬¥ç ­¨¥.} ¤­®§­ ç­ë¥ â®¯®«®£¨ç¥áª¨ § ¬ª­ãâë¥ ®â®¡à ¦¥­¨ï
ç áâ® ­ §ë¢ îâ ¯à®áâ® § ¬ª­ãâë¬¨.  â¥®à¨¨ ¬­®£®§­ ç­ëå
®â®¡à ¦¥­¨© § ¬ª­ãâë¬¨ ­ §ë¢ îâ ­¥áª®«ìª® ¤àã£¨¥ ®â®¡à ¦¥­¨ï (á¬.
¢ëè¥).

¨¬¢®« ¬¨ $R, N, Z$ ¡ã¤¥¬ ®¡®§­ ç âì á®®â¢¥âáâ¢¥­­® ¬­®¦¥áâ¢ 
¤¥©áâ¢¨â¥«ì­ëå, ­ âãà «ì­ëå, æ¥«ëå ç¨á¥«.

¢¥¤¥¬ ¥é¥ ¤¢  â¥à¬¨­ .

ë ¡ã¤¥¬ £®¢®à¨âì, çâ® ¬­®¦¥áâ¢® $A \subset X$ {\it à ¢­®¬¥à­®
â¥«¥á­®}, ¥á«¨ ¯¥à¥á¥ç¥­¨¥ $A$ á «î¡ë¬ ®âªàëâë¬ è à®¬ «¨¡® ¯ãáâ®,
«¨¡® á®¤¥à¦¨â ¤àã£®© ®âªàëâë© è à.

®­ïâ­®, çâ® «î¡®¥ ®âªàëâ®¥ ¬­®¦¥áâ¢® à ¢­®¬¥à­® â¥«¥á­®.
®«ã¨­â¥à¢ «ë ¨ ®âà¥§ª¨ ¢ $R$ â ª¦¥ à ¢­®¬¥à­® â¥«¥á­ë.

ë ¡ã¤¥¬ ­ §ë¢ âì ¬­®£®§­ ç­®¥ ®â®¡à ¦¥­¨¥ $F:X \multimap Y${\it
¯®çâ¨ ¯®«ã­¥¯à¥àë¢­ë¬ á­¨§ã}, ¥á«¨ ¤«ï «î¡®£® ®âªàëâ®£® $B \subset
Y$ ¬­®¦¥áâ¢® $ F_{ - }^{ - 1} (B)$ à ¢­®¬¥à­® â¥«¥á­®.

áïª®¥ ¯®«ã­¥¯à¥àë¢­®¥ á­¨§ã ¬ã«ìâ¨®â®¡à ¦¥­¨¥, ¢ ç áâ­®áâ¨,
®¤­®§­ ç­®¥ ­¥¯à¥àë¢­®¥ ®â®¡à ¦¥­¨¥, ¯®çâ¨ ¯®«ã­¥¯à¥àë¢­® á­¨§ã.

¬¥áâ¥ á ª ¦¤ë¬ $F:X \multimap Y$¬®¦­® à áá¬®âà¥âì ¬­®£®§­ ç­®¥
®â®¡à ¦¥­¨¥ $\widetilde{F} :F(X) \multimap X$, ª®â®à®¥
®¯à¥¤¥«ï¥âáï â ª: $\widetilde{F}(y) = F_{ - }^{ - 1} (y)$. ¥£ª®
¢¨¤¥âì, çâ® ¤«ï $B \subset Y:\widetilde{F}(B\bigcap F(X))=F_{ -
}^{ - 1} (B)$. â¬¥â¨¬ â ª¦¥, çâ® $\widetilde{F}_{ - }^{ - 1} (A)
= F(A)$ ¤«ï ¢á¥å $A \subset X$. ¥©áâ¢¨â¥«ì­®,

$$
y \in \widetilde {F}_{ - }^{ - 1} (A) \Leftrightarrow \widetilde
{F}(y) \cap A \ne \emptyset \Leftrightarrow $$ $$\Leftrightarrow
\exists x \in A, \, x \in F_{ - }^{ - 1} (y) \Leftrightarrow
\exists x \in A,F(x) \ni y \Leftrightarrow y \in F(A).
$$

\section{ à¨â¥à¨© á®¡áâ¢¥­­®áâ¨}

¥¯¥àì ¬ë áä®à¬ã«¨àã¥¬ ®á­®¢­®© à¥§ã«ìâ â à ¡®âë.

{\bf ¥®à¥¬  1.} {\it ãáâì X, Y -- ¬¥âà¨ç¥áª¨¥ ¯à®áâà ­áâ¢ , }
\\
$F:X \multimap Y${\it ¯®çâ¨ ¯®«ã­¥¯à¥àë¢­® á­¨§ã. ãáâì X ¯®«­®.
ãáâì ¤«ï ¢áïª®£® }$y \subset Y${\it ¬­®¦¥áâ¢® } $F_{ - }^{ - 1}
(y)$ {\it § ¬ª­ãâ®, ­® ­¥ á®¤¥à¦¨â ­¨ ®¤­®£® ®âªàëâ®£® è à .
®£¤ , ¥á«¨ F â®¯®«®£¨ç¥áª¨ § ¬ª­ãâ®, â® F á®¡áâ¢¥­­®.}

{\bf ®ª § â¥«ìáâ¢®.}

ãáâì $y_{0} \subset Y$ -- ¯à®¨§¢®«ì­ ï â®çª . ®ª ¦¥¬ á­ ç « ,
çâ® $F_{ - }^{ - 1} (y_{0} )$ ®â­®á¨â¥«ì­® ª®¬¯ ªâ­®. à¥¤¯®«®¦¨¬,
çâ® íâ® ­¥ â ª. ®£¤  ¯® â¥®à¥¬¥  ãá¤®àä  $F_{-}^{-1}(y_0)$ ­¥
¨¬¥¥â ª®­¥ç­®© $3\delta $ - á¥â¨, £¤¥ $\delta > 0$ ¤®áâ â®ç­®
¬ «®. ®íâ®¬ã $F_{ - }^{ - 1} (y_{0} )$ á®¤¥à¦¨â ¡¥áª®­¥ç­ãî
¯®á«¥¤®¢ â¥«ì­®áâì $z_{1} ,z_{2} ,...$ â®ç¥ª, ¯®¯ à­ë¥ à ááâ®ï­¨ï
¬¥¦¤ã ª®â®àë¬¨ ­¥ ¬¥­ìè¥ $3\delta$ (¥á«¨ ¡ë â ª¨å â®ç¥ª ¡ë«® «¨èì
ª®­¥ç­®¥ ç¨á«® $z_{1} ,...,z_{m} $, â® ®áâ «ì­ë¥ â®çª¨ $F_{ - }^{
- 1} (y_{0} )$ ®âáâ®ï«¨ ¡ë ®â íâ¨å ¬¥­ìè¥, ç¥¬ ­  $3\delta $; â.¥.
$z_{1} ,...,z_{m} $ ¡ë«® ¡ë $3\delta$-á¥âìî). ¬¥áâ¥ á ª ¦¤®©
$z_{i} (i \in N)$ à áá¬®âà¨¬ ®âªàëâë© è à $\widetilde {G}_{i} =
B(z_{i} ,\delta )$ á æ¥­âà®¬ ¢ $z_{i} $ à ¤¨ãá  $\delta $. á­®,
çâ® ¥á«¨ ¤¢¥ â®çª¨ «¥¦ â ¢ à §­ëå $\widetilde {G}_{i}$ ¨
$\widetilde {G}_{j} $, â® à ááâ®ï­¨¥ ¬¥¦¤ã ­¨¬¨ ­¥ ¬¥­¥¥ $\delta
$.  áá¬®âà¨¬ â¥¯¥àì ¢ $Y$ ®âªàëâë¥ è àë $B_{n} = B(y_{0}
,\frac{{1}}{{n}})$ $(n \in N) $. .ª. $F$ ¯®çâ¨ ¯®«ã­¥¯à¥àë¢­®
á­¨§ã, â® $F_{ - }^{ - 1} (B_{n} )$ à ¢­®¬¥à­® â¥«¥á­®.  áá¬®âà¨¬
¬­®¦¥áâ¢® $\widetilde{\widetilde{G}}_n=\widetilde{G}_n\bigcap
F_{-}^{-1}(B_n)$. ­® ­¥ ¯ãáâ®, â.ª. á®¤¥à¦¨â $z_{n} $. .ª. $F_{
- }^{ - 1} (B_{n} )$ à ¢­®¬¥à­® â¥«¥á­®, â® ­ ©¤¥âáï ª ª®©-­¨¡ã¤ì
®âªàëâë© è à, á®¤¥à¦ é¨©áï ¢ $\widetilde {\widetilde {G}}_{n} $.
¡®§­ ç¨¬ íâ®â è à ç¥à¥§ $G_{n} $. ® ãá«®¢¨î â¥®à¥¬ë $G_{n}$ ­¥
á®¤¥à¦¨âáï ¯®«­®áâìî ¢ $F_{ - }^{ - 1} (y_{0} )$, â.¥. ­ ©¤¥âáï
$x_{n} \in G_{n} ,y_{0} \notin F(x_{n} )$. ® $x_{n} \in G_{n}
\subset F_{ - }^{ - 1} (B_{n} )$. ®íâ®¬ã ­ ©¤¥âáï $y_{n} \in
B_{n} $ â ª®¥, çâ® $y_n\in F(x_n) $. .ª. â®çª¨ $x_{n} $ «¥¦ â ¢
$\widetilde {G}_{n} $, â® ¯®¯ à­ë¥ à ááâ®ï­¨ï ¬¥¦¤ã ­¨¬¨ ­¥ ¬¥­¥¥
$\delta $. .¥. $\{x_n\}_{n\in N}$ -- ¬­®¦¥áâ¢® ¡¥§ ¯à¥¤¥«ì­ëå
â®ç¥ª. ®íâ®¬ã ®­® § ¬ª­ãâ®. ® $F$ â®¯®«®£¨ç¥áª¨ § ¬ª­ãâ®,
¯®íâ®¬ã $F(\{ x_{n} \} _{n \in N}^{} )$ § ¬ª­ãâ®. ® $\{ y_{n} \}
_{n \in N}^{} \in F(\{ x_{n} \} _{n \in N}^{} )$. .ª. $y_{n} \in
B_{n} $, â® $y_{n} \mathop { \to }\limits_{n \to \infty } y_0 $;
¯®«ãç ¥âáï, çâ® $y_{0} \in F(\{ x_{n} \} _{n \in N}^{} )$. ®
íâ®£® ­¥ ¬®¦¥â ¡ëâì, â.ª. ¤«ï ¢áïª®£® $n$ $y_{0} \notin F(x_{n}
)$. à®â¨¢®à¥ç¨¥.

â ª, ¤«ï ¢áïª®£® $y_{0} $ ¬­®¦¥áâ¢® $F_{-}^{-1}(y_0)$
®â­®á¨â¥«ì­® ª®¬¯ ªâ­®. ® $F_{ - }^{ - 1} (y_{0} )$ § ¬ª­ãâ® ¯®
ãá«®¢¨î. ­ ç¨â, $\widetilde {F}(y_{0} ) = F_{ - }^{ - 1} (y_{0}
)$ ª®¬¯ ªâ­® ¤«ï ¢áïª®£® $y_{0} \in F(X)$.

 «¥¥, $\widetilde {F}:F(X) \multimap X$ ¯®«ã­¥¯à¥àë¢­® á¢¥àåã.
¥©áâ¢¨â¥«ì­®, ¥á«¨ $A$ § ¬ª­ãâ® ¢ $X$, â®
$\widetilde{F}_{-}^{-1}(A)=F(A)$ § ¬ª­ãâ® ¢ $Y$, â.ª. $F$
â®¯®«®£¨ç¥áª¨ § ¬ª­ãâ®. ç¥¢¨¤­® â®£¤ , çâ® $\widetilde {F}_{ -
}^{ - 1} (A)$ § ¬ª­ãâ® ¨ ¢ $F(X)$.

«¥¤®¢ â¥«ì­® (á¬. ¯.1.), ¯à¨ ¬ã«ìâ¨®â®¡à ¦¥­¨¨ $\widetilde {F}$
®¡à § «î¡®£® ª®¬¯ ªâ­®£® ¬­®¦¥áâ¢  ª®¬¯ ªâ¥­.

ãáâì $B$ -- ª®¬¯ ªâ­®¥ ¬­®¦¥áâ¢® ¢ $Y$. .ª. $F(X)$ § ¬ª­ãâ® ¢
$Y$, â® $F(X)\cap B$ ª®¬¯ ªâ­® (¢ $F(X)$). ®£¤  $\widetilde
{F}(F(X) \cap B)$ ª®¬¯ ªâ­®. .¥. $F_{ - }^{ - 1} (B)$ ª®¬¯ ªâ­®,
¨ á®¡áâ¢¥­­®áâì $F$ ¤®ª § ­ .

\section{ ¬¥ç ­¨ï}

{\bf  ¬¥ç ­¨¥ 1.}  ä®à¬ã«¨à®¢ª¥ â¥®à¥¬ë 1 ¯®çâ¨
¯®«ã­¥¯à¥àë¢­®áâì á­¨§ã ­¥«ì§ï § ¬¥­¨âì ­  ¯®«ã­¥¯à¥àë¢­®áâì
á¢¥àåã, çâ® ¯®ª §ë¢ ¥â á«¥¤ãîé¨© ¯à¨¬¥à.  áá¬®âà¨¬ ®â®¡à ¦¥­¨¥
$F:R \multimap R$:

 $F(x)=\{x\}$ (¬­®¦¥áâ¢®, á®áâ®ïé¥¥
¨§ ®¤­®© â®çª¨ $x$), ¥á«¨ $x \notin Z$;

$F(x) = R, $ ¥á«¨ $x \in Z$.

â®¡à ¦¥­¨¥ $F$ ¯®«ã­¥¯à¥àë¢­® á¢¥àåã ¨ â®¯®«®£¨ç¥áª¨ § ¬ª­ãâ®.
®«­ë© ¯à®®¡à § ¢áïª®© â®çª¨ $y$ ¥áâì ¬­®¦¥áâ¢® $\{ y\} \cup Z$.
®­ïâ­®, çâ® ®­® § ¬ª­ãâ® ¨ ­¥ á®¤¥à¦¨â ­¨ ®¤­®£® è à , ­® ­¥
ª®¬¯ ªâ­®. ®íâ®¬ã $F$ ­¥ á®¡áâ¢¥­­®.

{\bf  ¬¥ç ­¨¥ 2.} à¥¡®¢ ­¨¥ ¯®çâ¨ ¯®«ã­¥¯à¥àë¢­®áâ¨ á­¨§ã ¢
â¥®à¥¬¥ 1 áãé¥áâ¢¥­­® á« ¡¥¥ âà¥¡®¢ ­¨ï ¯®«ã­¥¯à¥àë¢­®áâ¨ á­¨§ã.
¥©áâ¢¨â¥«ì­®, à áá¬®âà¨¬ á«¥¤ãîé¨© ¯à¨¬¥à. ãáâì ®â®¡à ¦¥­¨¥ $F:R
\multimap R$ § ¤ ­® â ª:

 $F(x) = \{ \left| {x} \right| + 1\} $ ¯à¨ $\left| {x} \right| > 1$;

  $F(x) = \{ \left| {x} \right| + 1,1 - \left| {x} \right|\} $ ¯à¨ $\left| {x}
\right| = 1$;

 $F(x)=\{ \left| {x} \right| + 1,1 - \left| {x}
\right|,\left|{x}\right|-1\}$ ¯à¨ $\left| {x} \right| < 1$.

à¨ ®â®¡à ¦¥­¨¨ $F$ ¯®«­ë© ¯à®®¡à § ¢áïª®© â®çª¨ ª®­¥ç¥­,   ¯®â®¬ã
§ ¬ª­ãâ ¨ ­¥ á®¤¥à¦¨â è à . á­® â ª¦¥, çâ® $F$ â®¯®«®£¨ç¥áª¨
§ ¬ª­ãâ®. ®«­ë© ¯à®®¡à § «î¡®£® ®âªàëâ®£® ¬­®¦¥áâ¢  à ¢­®¬¥à­®
â¥«¥á¥­, ­® ­¥ ®¡ï§ â¥«ì­® ®âªàëâ: ­ ¯à¨¬¥à, ¯®«­ë© ¯à®®¡à § ¬ «®©
®ªà¥áâ­®áâ¨ ­ã«ï á®¤¥à¦¨â ªà ©­¨¥ â®çª¨ +1 ¨ -1. â ª, ®â®¡à ¦¥­¨¥
$F$ ã¤®¢«¥â¢®àï¥â ãá«®¢¨ï¬ â¥®à¥¬ë, ­® ­¥ ï¢«ï¥âáï ¯®«ã­¥¯à¥àë¢­ë¬
á­¨§ã.

\section{¡à â­ ï â¥®à¥¬ }

â®â ¯ã­ªâ ¯®á¢ïé¥­ â¥®à¥¬¥, ª®â®àãî ¬®¦­® ­ §¢ âì ®¡à â­®© ª â¥®à¥¬¥ 1.

{\bf ¥®à¥¬  2. }{\it ãáâì X, Y -- ¬¥âà¨ç¥áª¨¥ ¯à®áâà ­áâ¢ , }
\\
$F:X \multimap Y${\it -- § ¬ª­ãâ®¥ ¬ã«ìâ¨®â®¡à ¦¥­¨¥. ®£¤ , ¥á«¨
F á®¡áâ¢¥­­®, â® F â®¯®«®£¨ç¥áª¨ § ¬ª­ãâ®.}

{\bf ®ª § â¥«ìáâ¢®.}

ãáâì $A$ - § ¬ª­ãâ®¥ ¬­®¦¥áâ¢® ¢ $X$. ®ª ¦¥¬, çâ® $F(A)$
§ ¬ª­ãâ®. ¥©áâ¢¨â¥«ì­®, ¯ãáâì $y_{n} \in F(A)$, $y_{n} \mathop {
\to }\limits_{n \to \infty } y $. ­®¦¥áâ¢® $Y_{0} = \{ y\} \cup
\{y_{n}\}_{n \in N} $ ª®¬¯ ªâ­®. ®íâ®¬ã ¨ $F_{ - }^{ - 1} (Y_{0}
)$ ª®¬¯ ªâ­®. .ª. $y_{n} \in F(A)$, â® ­ ©¤¥âáï $x_{n} \in A$
â ª®¥, çâ® $y_n\in F(x_n)$. ®­ïâ­®, çâ® $\{ x_{n} \} _{n \in
N}^{} \in F_{ - }^{ - 1} (Y_{0} )$. ®íâ®¬ã $\{ x_{n} \} _{n \in
N}^{} $ ¡¥§ ®£à ­¨ç¥­¨ï ®¡é­®áâ¨ áå®¤¨âáï; $x_{n} \mathop { \to
}\limits_{n \to \infty } x$. .ª. $x_{n} \in A$, â® $x\in A$. ®
$F$ § ¬ª­ãâ®, ¯®íâ®¬ã $y \in F(x) \subset F(A)$. â ª, $F(A)$
§ ¬ª­ãâ®. ¥®à¥¬  ¤®ª § ­ .

\bigskip

{\bf ¨â¥à âãà .}

\bigskip

1. ®à¨á®¢¨ç .., ¢ï£¨­ ..,  ¯à®­®¢ .. ¥«¨­¥©­ë¥ äà¥¤£®«ì¬®¢ë
®â®¡à ¦¥­¨ï ¨ â¥®à¨ï ¥à¥ -  ã¤¥à  // , 1977. - .32. - ë¯.4. - .3-54.

2. ®à¨á®¢¨ç .., ¥«ì¬ ­ .., ëèª¨á .., ¡ãå®¢áª¨© .. ¢¥¤¥­¨¥ ¢
â¥®à¨î ¬­®£®\-§­ ç­ëå ®â®¡à ¦¥­¨©. - ®à®­¥¦: §¤. ®à®­¥¦áª®£® ã­¨¢¥àá¨â¥â ,
1986.-104á.

\end{document}